\documentclass{amsart}
\usepackage[mathscr]{eucal}

\usepackage{amsfonts}
\usepackage{amsmath}
\usepackage{amsthm}
\usepackage{amssymb}
\usepackage{latexsym}

\newtheorem{theorem}{Theorem}[section]
\newtheorem{lemma}[theorem]{Lemma}

\newtheorem{proposition}[theorem]{Proposition}
\newtheorem{corollary}[theorem]{Corollary}
\theoremstyle{definition}
\newtheorem{definition}[theorem]{Definition}

\def\leq{\leqslant}
\def\geq{\geqslant}
\def\eps{\varepsilon}
\def\epso{\varepsilon_1}

\def\T{\mathbb T}

\def\dM{\partial M}
\def\dMh{\partial\hat{M}}

\def\fh{\widehat{F}}
\def\dOp{\partial\Omega^+}
\def\R{\mathbb R}
\def\RR{\mathcal R}

\def\tang{\mathcal{S}}
\def\s3{\mathbf S^3}

\def\T{\mathbb T}
\def\sing{T^{-1} \mathcal{S}}
\def\ttang{T^{-2}\mathcal{S}}

\def\U{\mathcal U}

\def\Z{\mathbb Z}

\def\TM{\mathrm{T}_1\R^3}
\newcommand{\da}{\partial^{\alpha}}
\def\geq{\geqslant}
\def\leq{\leqslant}

\def\phi{\varphi}

\def\G{\mathcal{G}^\delta}
\def\Gi{G_i^{\delta}}

\def\be{\begin{equation}}
\def\ee{\end{equation}}

\def\bea{\begin{eqnarray}}
\def\eea{\end{eqnarray}}

\def\ba{\begin{array}}
\def\ea{\end{array}}

\begin{document}

\title{On the structure of the singularity manifolds of dispersing billiards}

\author{Pavel Batchourine}
\address{Department of Mathematics, Princeton University,
Princeton, NJ 08544, USA} \email{pbatchou@math.princeton.edu}

\maketitle

\renewcommand{\theequation}{\arabic{section}.\arabic{equation}}

\section{Introduction.}
A billiard is a dynamical system describing the motion of a point particle in a
connected, compact domain $Q\subset\T^d,$ $d\geq 2.$ In general, the boundary
of the domain $Q$ is assumed to be piecewise $C^3$-smooth, however later on we
impose further restrictions on the boundary (cf. Section $2$). Smooth
components of the boundary $\partial Q$ are called {\it scatterers}. If the
smooth components of the boundary are strictly convex, the billiard is called
{\it dispersing} or a {\it Sinai billiard}. Inside $Q$ the particle moves with
constant velocity with elastic reflections at the boundary. The kinetic energy
of a point is a first integral of the motion. Therefore the phase space of the
billiard is $M=Q\times S^{d-1}=\{(q,v)|q\in Q,\ ||v||=1\}.$ The Liouville
probability measure $\mu$ on $M$ is a product of the Lebesgue measures on $Q$
and $S^{d-1},$ i.e. $d\mu=const\cdot dq dv.$ The resulting dynamical system
$(M,\{S^t,t\in\R\},\mu)$ is the billiard flow.

The boundary $\partial Q$ defines a natural cross-section for the billiad flow.
Let $n(q)$ be the unit normal vector of a smooth component of the boundary
$\partial Q$ at the point $q\in\partial Q$, directed inwards. The {\it billiard
map} $T$ is defined as the first return map on $\partial M,$ where
$$\partial M=\{(q,v)|q\in\partial Q,\ (v,n(q))\geq 0\}.$$

The invariant measure for the map $T$ is denoted by $\mu_1,$ and we have
$$d\mu_1=const\cdot |(v,n(q))|dq dv.$$

Contrary to the case of smooth dynamical systems, the billiard map $T$ is
discontinuous, which make the application of classical methods of analysis of
stochastic properties of a system significantly more difficult. Nevertheless,
in $1970$ Sinai demonstrated (\cite{S1}) that at least in the case $d=2$
dispersing billiards are hyperbolic, ergodic and even $K$-mixing. The main tool
in the proof of ergodicity of the billiard system was the so-called {\it
fundamental theorem} for the dispersing billiards. This theorem ensures an
abudance of geometrically nicely situated and sufficiently large stable and
unstable manifolds. For the case $d>2$ the fundamental theorem was first stated
and proved in \cite{SCh} and generalized in \cite{KSSz}.

In the proof of the fundamental theorem for dispersing billiards one makes some
assumptions concerning the structure of the set of points, where the billiard
map and its iterations are discontinuous. These sets are called singularity
manifolds. They were assumed to be smooth manifolds with boundary until
recently, when it was found in \cite{BCST1} that in a typical situation there
are special points where these manifolds may have singularities in the sense of
the singularity theory. At these points the singularity manifolds don't have a
tangent plane. Therefore one has to analyze the singularities to show that the
arguments in the proof of the fundamental theorem remain valid. It was shown in
\cite{BCST2} that in the case of algebraic semi-dispersing billiards the proof
of the fundamental theorem remains valid.

In this paper we show that in the case of strictly dispersing scatterers
in $Q\subset\T^d,$ $d\geq 3$ the
assumption about algebraicity is not needed and the fundamental theorem remains
valid if the boundaries of the scatterers are $C^{8d-7}$ smooth.

 The paper consists of two parts. First, we analyze the singularity manifolds
of the billiard map. These manifolds can be obtained as the images of some
smooth manifold $M$ under several iterations of the billiard map. This map is
irregular at the the points of tangency, so to understand the structure of the
singularity one needs to study the billiard map near the points of tangency. We
compute partial derivatives of the map and show (Lemma~\ref{change}), that
after a singular, but simple change of variables the billiard map becomes a
local diffeomorphism. In Lemma \ref{resolution} we show that the singularity
manifold is a level set of a smooth function with a non-vanishing germ. The
Theorem \ref{measure} implies that the singularity manifold has a certain
regularity property required in the proof of the fundamental theorem.

In the second part we show that the proof of \cite{SCh} of the fundamental
theorem remains valid for typical configurations of strictly dispersing
scatterers. The important observation here is Lemma \ref{4ref} which states
that for a typical configuration strictly dispersing scatterers the
trajectories can have at most $2d-2$ tangencies.

\subsection*{Acknowledgements}

I would like to thank Ya.G. Sinai for introducing me to this area and proposing
the original problem. I am very grateful to him for his constant interest,
encouragement and help.

\setcounter{equation}{0}
\section{Preliminaries}

\subsection{Assumptions about the scatterers}
We make the following assumptions about the the scatterers:

\begin{enumerate}
    \item There are finitely many scatterers $\RR_1,\ldots,\RR_K$,
    \item Each scatterer is a strictly convex $C^{8d-7}$
    hypersurface,\label{convex}
    \item Finite horizon: the length of the trajectory between two successive
        reflections is uniformly bounded,
    \item Time between two successive reflections is bounded from below by some
            positive constant $\tau_0>0.$\label{nomult}
\end{enumerate}

Strictly convex here means that the second fundamental form of each scatterer
is positively defined at every point.

Throughout the paper we consider the configuration space $Q\subset\T^d$ as a
subset of the $\R^d,$ the universal cover of $\T^d.$ The configuration of
scatterers in $\T^d$ lifts to a periodic configuration of scatterers in $\R^d.$
We denote the phase space of the corresponding billiard system in $\R^d$ by
$\dMh$ and use the same notation for the billiard map in $\R^d:$
$T:\dMh\to\dMh.$

To fix the notations for the scatterers we choose $K$ scatterers
$\widehat{\RR}_1,\ldots \widehat{\RR}_K$ in $\R^d,$ so that these
scatterers
satisfy the following two properties:

\begin{itemize}
\item for every scatterer $\RR_m$ on $\T^d$ there exists a unique scatterer
    $\widehat{\RR}_m$ in $\R^d,$ which projects to $\RR_m$
\item $\widehat{\RR}_m$ is the closest to the origin among the preimages of $\RR_m$
\end{itemize}

Fix coordinate system in $\R^d$ and assume that $\RR_m$  is given as zero level
set of a function $R_m\in C^{8d-7}(\R^d).$ The rest of scatterers are then
given as zero level sets of the shifts of functions $R_k$ by elements of the
lattice which generates the torus. The assumption (\ref{convex}) on the
scatterers implies that for each $k=1,\ldots K$ the hessian of $R_k$ is
positive definite at each $z\in\R^d,$ such that $R_k(z)=0.$

\subsection{Singularities} Consider the set of tangential reflections
$$\mathcal{S}:=\{(q,v)\in\partial M | (v,n(q))=0\}.$$

It is easy to see, that the map $T$ is not continuous at the set $\sing.$ As a
result, the singularity set for any iterate $T^n$ is
$$\mathcal{S}^{(n)}=\cup_{i=1}^{n}\mathcal{S}^{-i},$$
where in general $\mathcal{S}^k=T^k\mathcal{S}.$ The set of regular
trajectories, i.e. the trajectories, for which $T^i x\notin \mathcal{S},$
$-\infty<i<\infty,$ will be denoted by
$$\partial M^0:=\partial M\diagdown \cup_{n\in\Z}\mathcal{S}^n.$$

In the case when the boundary $\partial Q$ is only piecewise smooth, additional
singularities - multiple collisions - arise. They correspond to points
$q\in\partial Q,$ which belong to several smooth components of the boundary. At
these points the dynamics is not well-defined. The assumption (\ref{nomult})
implies that we consider only {\it tangential} singularities. It is enough to
consider this situation, since the blow-up of the derivative of the map $T$ -
the effect found in \cite{BCST1} doesn't occur for the singularities
corresponding to multiple collisions.

\section{Space of lines} The dynamics in a billiard system consists
of two parts: the motion along a straight line and a reflection. It is the
reflection part that causes discontinuities and singularities, so it is natural
to isolate it. One way to do so is to look at the dynamics on the space of
oriented lines.

In this and next section we assume, that the dynamics takes place in $\R^d.$
Consider a point $x=(q,v)\in\partial\hat{M}$ in the phase space of a billiard.
The pair $(q,v)$ uniquely determines an oriented line $l$ in $\R^d.$ The image
$x'=(q',v')=T(x)$ determines a new line $l'$, which is by definition the image
of $l.$ This map on the space of oriented lines will be called the billiard map
on the space of oriented lines. Since the same oriented line can intersect
several scatterers, the dynamics is globally not well-defined. However, if we
are interested in what happens near a fixed trajectory in a finite number of
iterations of the billiard map, then we know from which scatterers the lines
are reflected and so the dynamics is well-defined.

\subsection{Coordinates on the space of oriented lines}
A line $l$ in $\R^d$ can be characterized by its closest point to the origin
$p$, and by its direction $v.$ Vector $v$ can be normalized to have length
one, since multiplying $v$ by any positive constant still gives the same
direction of the line: \be ||v||=1\label{cond1}\ee Clearly,
\be(p,v)=0\label{cond2}\ee These coordinates give the identification of he
space of oriented lines in $\R^d$ with the (co)tangent bundle of unit sphere
$S^d.$ We denote this space by $\Omega.$ The symplectic structure given by this
identification differs from the standard symplectic structure on $T^* S^d.$ by
sign. For a fixed $v$ we shall denote the plane (\ref{cond2}) in $(\R^d,p)$ by
$\Pi_v.$

\begin{lemma}
Suppose $x_0=(q_0,v_0)\in\dM$ and $U(x_0)$  is a neighborhood of $x_0$ in
$\dM,$ such that for every $x=(q,v)\in U(x_0)$ we have

\be (v,n(q))\geq c > 0 \label{away}, \ee

where $c$ is some positive constant. Then there exists a diffeomorphism
$\Phi\in C^{8d-7}(U(x_0))$
from the neighborhood $U(x_0)$ to the neighborhood of the line $l_0$ determined
by $x_0,$ such that $\Phi(x)=l_0$ and $|det D\Phi|\geq c.$ \label{tolines}
\end{lemma}
\begin{proof}
The diffeomorphism $\Phi$ is given by the formula
$$\Phi(q,v)=(\pi_v(q),v),$$
where $\pi_v:\R^d\to\R^d$ is a projection to the plane $\Pi_v$ in the direction
of vector $v.$ It follows from (\ref{away}), that $|det D\Phi|\geq c.$
\end{proof}

\section{Local structure of a reflection near tangency}
In this section we consider the reflection of trajectories from a fixed
scatterer given by the equation $F(x_1,\ldots,x_d)=0.$ We start by showing that
the set of lines tangent to this scatterer is a smooth hypersurface in the
space of oriented lines.

Consider a unit tangent bundle $\mathrm{T}_1\R^d.$ It is obviously a bundle
over $\R^d$ with fibers isomorphic to $S^d.$ Any function $F(x_1,\ldots,x_d)$
lifts to a function $\fh$ on $\mathrm{T}\R^d,$ which constant along the fibers.
The space $\mathrm{T}_1\R^3$ is also a bundled over $\Omega,$ the space of
lines in $\R^d.$ Denote by $\pi$ the projection corresponding to this bundle.
The fibers $\pi^{-1}(\omega),$ $\omega\in\Omega$ are the points in $\R^d$ which
lie on a given line with the direction determined by $\omega.$

Denote the coordinate in the fiber by $t.$ If $x=(p,v)$ are coordinates in
$\Omega,$ then $(x,t)$ are the coordinates in $\mathrm{T}_1\R^d.$ The set of
lines tangent to the scatterer is the projection of the discriminant set
$\Sigma$ given by $$\fh(x,t)=0,\quad \frac{\partial \fh}{\partial t}(x,t)=0$$
on $\Omega.$

\begin{lemma}
Let $(x_0,t_0)\in\Sigma.$ There exists a neighborhood $U$ of $(t_0,x_0)$ in
$\R\times\Omega,$ such that $\pi(\Sigma\cap U)$ is a parameterized $2d-3$
dimensional $C^{8d-9}$ manifold in $\Omega.$
\label{smoothtan}
\end{lemma}

\begin{proof}

Consider a mapping $G:\R^{2d-2}\times\R\to\R^2,$ given by $G(x,t)=(\fh,
\frac{\partial \fh}{\partial t})(x,t).$ Then $\Sigma=G^{-1}(0).$ Jacobi matrix
for $G$ is equal to
\[\left( \begin{array}{cccc}
\frac{\partial \fh}{\partial x_1} & \ldots & \frac{\partial \fh}{\partial x_{2d-2}} &
\frac{\partial \fh}{\partial t}\\
\frac{\partial \fh}{\partial x_1 \partial t} & \ldots & \frac{\partial
\fh}{\partial x_{2d-2}\partial t} & \frac{\partial^2 \fh}{\partial t^2}
\end{array}\right)\]

Notice that from the assumption (\ref{convex}) on the scatterers, it follows
that $$\partial^2 \fh / \partial t^2 (x_0,t_0) \neq 0$$ Also, there exists an
index $i,$ $1\leq i\leq d-1,$ such that $\partial\fh /{\partial
x_i}(x_0,t_0)\neq 0.$ Assume that $i=d-1.$ Then by the Implicit Function
Theorem applied to $G$ we can express $t$ and $x_{d-1}$ as the functions of
$x_1,\ldots,\widehat{x_{d-1}},\ldots, x_{2d-3}$ on $\Sigma\cap U$ for a
neighborhood $U$ of $(t_0,x_0).$ Hence, on $\pi(\Sigma\cap U),$ $x_{d-1}$ is a
smooth function of $x_1,\ldots,\widehat{x_{d-1}},\ldots, x_{2d-2}.$
\end{proof}

Function $\fh:\Omega\times\R\to\R$ induces a function $\bar{F}:\Omega\to\R.$
For a given $\omega\in\Omega$
$$\bar{F}(\omega)=\inf\limits_{(x,t)\in\pi^{-1}(\omega)}\fh(x,t).$$

It follows from the assumption (\ref{convex}) on the scatterers that for
$\omega$ sufficiently close to $\dOp$ the infinum is attained at a unique point
$(\omega,t(\omega).$ The zero level set of $\bar{F}$ coincides with the set of
lines tangent to the scatterer. Let
$\Omega^+=\{\omega\in\Omega|\bar{F}(\omega)\leq0\}$ be the set of lines that
experience a reflection from the scatterer. Its boundary $\dOp$ consists of
lines tangent to the scatterer. On $\Omega^+$ the billiard map $T$ is defined.
It sends the lines before the reflection from the scatterer to the reflected
lines.

Let $x_0\in\dOp.$ By the previous Lemma there exists a neighborhood
$U(x_0)\subset\Omega,$ such that the $\dOp\cap U(x_0)$ is a hypersurface
parameterized by $(x_1,\ldots,\widehat{x_{d-1}},\ldots x_{2d-2}).$ In
particular, this implies that for any fixed $v_1$ sufficiently close to $v_0$
the set of lines $\{(p,v)\in U(x_0)\cap\dOp\ |\ v=v_1\}$ is a parameterized
$d-2$ dimensional manifold in $\Pi_{v_1}$. We can assume that it is
parameterized by $(p_1,\ldots p_{d-2}).$ For every point on this manifold we
consider a line passing through this point in the direction perpendicular to
$\dOp\cap U(x_0).$ Let's choose a parameter $\tau$ on these lines so that
$\tau=0$ corresponds to a point on $\dOp$ and so that $\tau>0$ corresponds to
reflection. For any $x\in U(x_0)$ there exists a point $x_1\in\dOp\cap U(x_0)$
and $\tau,$ such that $x$ lies on the line passing through $x_1$ perpendicular
to $\dOp$ and $distance (x_1,\dOp\cap U(x_0))=\tau.$ Let
$$\lambda_0=\min\limits_{x\in\R^d: F(x)=0}\min\limits_{w\in\R^d:||w||=1}(Hess(F)(x)w,w),$$
where $Hess(F)(x)$ is the Hessian of $F$ at $x.$ Pick a neighborhood
$U_1(x_0)\subset U_0(x_0)$ such that for every $x\in U_1(x_0)$ we have
$distance(x,\dOp)<\lambda_0/2.$ Then the neighborhood $U_1(x_0),$ is
parameterized by $(p_1,\ldots p_{d-2},\tau,v).$

\begin{definition} Suppose $H^+=\{(x_0,\ldots x_n)|x_n>0\}.$ The mapping
\be\left\{ \begin{array}{l}
y_0=x_0, \\
\ \ \ \ldots \\
y_{n-1}=x_{n-1} \\
y_n=\sqrt{x_n}
\end{array}\right.\label{qrchange}\ee
of $H^+$ into itself will be called a {\it quasi-regular} change of variables.
\end{definition}

\begin{lemma}
    After a quasi-regular change of variables $\upsilon=\sqrt{\tau}$ in a semi-neighborhood
    $U_1(x_0)\cap\Omega^+$ parameterized as above the billiard map becomes a $C^{4d-4}$
    diffeomorphism.
    \label{change}
\end{lemma}

\begin{proof}

The lemma follows from the following decomposition of the billiard map.

Consider again the unit tangent bundle $\TM$ as a bundle over the space of
lines $\Omega.$ The scatterer $F=0$ defines a section of this bundle over the
points corresponding to lines that experience a reflection: to each such lines
one associates the point on the scatterer where the line hits the scatterer and
the direction of this line. The points of this section are the pairs
$\{(p,v)|p\in\R^d,\ R(p)=0,\ ||v||=1\}$ form a surface $\mathcal{C}$
topologically equivalent to $S^d\times S^d.$ On this surface the following
smooth map is defined. To each $(p,v)\in\mathcal{C}$ one associates $(p',v'),$
where $p=p'$ and $v'$ is a reflection of $v$ with respect to $\nabla R(p).$
This map is smooth on the surface $\mathcal{C},$ but the projection from
$\Omega$ to $\mathcal{C}$ is not smooth at the points corresponding to the
lines tangent to $R=0.$ The only thing one needs to show is that this
projection is smooth after a quasi-regular change of variable.

Fix any direction $v\in S^d$ close to $v_0.$ The projection from $\Omega$ to
$\mathcal{C}$ for a fixed $v$ is the projection from the hyperplane $\Pi_v$ to
the scatterer in $\R^d.$ Let $t$ be the coordinate in $\R^d$ in the direction
perpendicular to $\Pi_v.$ For fixed coordinates $p_1,\ldots,p_{d-2}$ the
projection is a function $t=f(\tau).$ Assume that $f(0)=0.$ Since the scatterer
is strictly convex at each point, we get that $\tau=t^2 g(t),$ where $g(0)\neq
0$ and so $\upsilon=\pm t\sqrt{g(t)}.$ We choose the sign so that it
corresponds to the trajectories coming to the scatterer. Then by the Implicit
Function  Theorem $t$ is a smooth function of $\upsilon$ and so after a
quasi-regular change of coordinated the billiard map becomes smooth.
\end{proof}

\section{Resolution of the singularity}

Using the description of the billiard map given by Lemma (\ref{change}) we show
that the preimage of a level set of some $F\in C^k(\Omega)$ function under the
billiard map can be represented as a level set of a function of the same class
of smoothness.

The billiard map (and it's inverse) is smooth away from tangencies, so in a
neighborhood of any point corresponding to a non-tangential reflection the
pull-back of $F$ is smooth. We show that in a (semi)neighborhood of a tangency
the preimage of a level set of a smooth function can still be given as a level
set of a smooth function. Consider the coordinates from the Lemma~\ref{change}
and suppose that the set of lines (the example in mind is the set of lines
tangent to some other scatterer or its preimage) is given by the equation
$$F(x_1,\ldots x_{2d-2})=0$$

In the the decomposition given by the lemma, only the substitution
$\tau\to\sqrt{\tau}$ is singular, so it is enough to look at this part
separately.

\begin{lemma}
    Let $\U$ be a neighborhood of zero in $\R^n,$
    $F\in C^k(\U\cap\{(x_1,\ldots,x_n)|x_n\geq 0\}),$ then there
    exists a neighborhood $\U_1\subset\U$ of zero and functions
    $G_{+},$ $G_{-}\in C^k(\U_1),$ such that
    \be
        F(x_1,\ldots,x_n)=
            G_{+}(x_1,\ldots,x_n^2)+x_n G_{-}(x_1,\ldots,x_n^2)
    \label{decomp}
    \ee
\label{resolution}
\end{lemma}

\begin{proof}
Extend $F$ (for example, by Whitney extension theorem) to small full
neighborhood of zero and denote the continuation by $\hat{F}.$ Let
$$F_{+}(x)=\frac{\hat{F}(x)+\hat{F}(-x)}{2},\quad F^{-}=
\frac{\hat{F}(x)-\hat{F}(-x)}{2}$$ be the decomposition of $\hat{F}$ in even
and odd parts. By a theorem of Glaeser \cite{BL} a function, which is even in
some variable is a function of the square of this variable, so there
exists a neighborhood $\U_1,$ and functions $G_{+},$ $G_{-}$ as in the
statement of the lemma, such that $F_{+}(x_1,\ldots x_n)=G_{+}(x_1,\ldots,
x_n^2)$ and $F_{-}(x_1,\ldots, x_n)=x_n G_{-}(x_1,\ldots, x_n^2).$
This implies the existance of the decomposition (\ref{decomp}).
\end{proof}

\begin{corollary}
 The level set $\{F(x_1,\ldots,\sqrt{x_n})=0\}\cap\U_1$ is given by
the equation $G_{+}^2=x_n G_{-}^2.$
\end{corollary}

It follows from (\ref{decomp}) that if $F$ has a non-zero $m$-jet,
$0<m<k/2$ then $G_{+}^2-x_n G_{-}^2$ has non-zero $2m$-jet.

\section{No more than $2d-2$ tangencies for generic scatterers}
Since the dimension of the phase space is equal to $2d-2,$ it is natural to
assume that, in a certain sense, no degenerations of codimension greater than
$2d-2$ can appear in a typical billiard system. Here we formulate and prove one
such statement that will be used in the next section.

Recall that the billiard system in $\R^d$ consists of $K$ scatterers given by
the equation $R_i=0,$ $i=1,\ldots K$ and their shifts by the elements of the
lattice which generates the torus $\T^d.$ The collection of functions
$\{R_1,\ldots R_K\}$ is a point in $[C^{4d-4}(\R^d)]^K.$

\begin{proposition}
    For any $0\leq m\leq 8d-7$ and $\eps>0$ there exists a collection of
scatterers
    $\{R_1(\eps),\ldots R_K(\eps)\}\in[C^{8d-7}(\R^d)]^K,$ such that
    \begin{itemize}
        \item for the billiard system defined by $\{R_1(\eps),\ldots R_K(\eps)\}$ there are
            no trajectories with more than $2d-2$ (not necessarily successive) tangencies
        \item $||R_1(\eps)-R_1||^2_{C^m(\R^d)}+\ldots+||R_K(\eps)-R_K||^2_{C^m(\R^d)}<\eps$
    \end{itemize}
\label{4ref}
\end{proposition}
\begin{proof}

We prove the proposition for the case $d=3$ and some fixed $m\geq 2.$  The
proof proceeds by induction on the number of reflections that can occur between
tangencies. At the step $k$ we make perturbations of the scatterer, so that
\begin{itemize}
\item The set of trajectories with three tangencies and at most $k$ proper
reflections in between form a finite union of ruled surfaces in the
configuration space.
\item Each ruled surface intersects every scatterer transversally
\item There are no $5$-tangencies with at most $k$ reflections in between
\end{itemize}

Let's first describe two types of perturbations used in the proof.
\begin{enumerate}
\item[Type 1]
Suppose that a ruled surface in $\R^3$ is tangent to a smooth strictly convex
surface given by $F(x,y,z)=0$ along a curve $\gamma.$ Then for every
$\eps_1,\eps_2>0$ there exists a smooth strictly convex function $\tilde{F},$
such that $||F-\tilde{F}||_{C^m(\R^3)}<\eps_1,$ $F(x,y,z)=\tilde{F},$ if the
distance from $(x,y,z)$ to $\gamma$ in $\R^3$ is greater or equal to $\eps_2$
and such that the ruled surface intersects the surface $\tilde{F}=0$
transversally at each point. In a degenerate situation, when $\gamma$ is a
point, this type includes the case when a ruled surface has a non-transversal
intersection with the convex surface at one point.

\item[Type 2]
Suppose that a line $l$ is tangent to a smooth strictly convex surface given by
$F(x,y,z)=0$ at a point $p.$ Then for every $\eps_1,\eps_2>0$ there exists a
smooth strictly convex function $\tilde{F},$ such that
$||F-\tilde{F}||_{C^m(\R^3)}<\eps_1,$ $F(x,y,z)=\tilde{F},$ if the distance
from $(x,y,z)$ to $p$ in $\R^3$ is greater or equal to $\eps_2$ and such that
the line $l$ intersects the surface $\tilde{F}=0$ transversally.
\end{enumerate}



Let us call by a combinatorial type of a trajectory of length $L$ a sequence of
$L$ scatterers in $\R^n$ labelled by "what happens" at this scatterer (a
reflection or a tangency).

We know from Lemma (\ref{smoothtan}), that the set of trajectories tangent to a
given scatterer form a smooth hypersurface. In the same way, the set of
trajectories which have two successive tangencies to given two scatterers is a
smooth $2$-dimensional manifold $\mathcal{T}$ with a boundary. If we allow
additional reflections between the two tangencies, then additional components
appear. Due to the effect described in \cite{BCST1} the manifold may behave
irregular at certain points of new components of the boundary. However, the set
of these points $\mathcal{S}$ is a finite union of smooth one-dimensional
manifolds. Hence by small perturbation of type $1$ of the scatterers we can
make the ruled surfaces which correspond to these manifolds transversal to the
scatterers and so that no points of $\mathcal{S}$ belong to more than one
component of the boundary.

This implies, that if we consider the set of trajectories tangent to the
previous two scatterers (with a number of reflections in between) and a new
scatterer, then (possible after an additional perturbation of type $1$ of the
new scatterer) the set of tangents to the new scatterer intersect $\mathcal{T}$
along a curve which has no points at which $\mathcal{T}$ is irregular.

This shows that the set of trajectories with three tangencies of any
combinatorial type of fixed length is a one-parameter family of trajectories
which form a ruled surface in the configuration space. The set of trajectories
with $4$ tangencies of a fixed combinatorial type maybe non-discrete if and
only if this ruled surface (or it's image under a billiard map) is tangent to
some scatterer along some curve on the scatterer. We make the perturbations of
Type $2$ of the scatterers so that all such ruled surfaces are transversal to
all scatterers using perturbations of type $1.$

Thus for every fixed length of combinatorial type we get finitely many
trajectories with $4$ tangencies. Suppose there exists a $5$-tangency among
them. Consider a line which corresponds to the trajectory "right before" the
fifth tangency. By choosing a sufficiently small perturbation of the scatterer
at the point of the fifth tangency we can make the line transversal to the
scatterer and not create any new $5$-tangencies. Thus we can eliminate
$5$-tangencies for a fixed length of the combinatorial class one by one.

To make this work we use the induction on the number of reflections between the
first and the last tangency.

Let $k$ denote the length of the considered trajectories. We start with $k=3.$
Fix some $\eps_0>0$ First, we find some ruled surface of order $0,$ which is
tangent to some scatterer along a curve. There is a finite number of ruled
surfaces of order $0$ associated with this scatterer. Among them some of the
ruled surfaces are transversal to each scatterer at each point of intersection.
We find a perturbation which destroys this non-discrete tangency and keeps all
ruled surface of order $0,$ which are associated to this scatterer transversal
to other scatterers. Now the number of ruled surfaces of order $0$ which are
tangent to some scatterer along a curve decreased by one. In the same way we
eliminate all non-transversal intersections of ruled surfaces with the
scatterers at points. This way we get that the set of trajectories with four
successive tangencies is finite.

Notice that there exists a positive number $\alpha,$ such that any perturbation
of size less that $\alpha,$ all ruled surfaces for the perturbed system are
transversal to every scatterer. Here $\alpha$  depends on the angles between
ruled surfaces of order zero and the scatterers.

We now make additional perturbations to eliminate all five successive
tangencies. Fix one successive $5$-tangency. Since the set of successive
$4$-tangencies is finite, there exists a neighborhood $\mathcal{U}$ of the
point of tangency on the last scatterer, in which the minimal angle between the
directions of incoming $4$-tangencies and the trajectories tangent to the
scatterer at the points in $\mathcal{U}$ is greater or equal to $c_1.$ Then by
perturbing (Type $2$) the scatterer in $\mathcal{U}$  so that the directions of
the lines tangent to the scatterer at the points of $\mathcal{U}$ change by at
most $c_1/2$ we make the scatterer transversal (locally) to the trajectory and
we don't create additional successive $5$-tangencies are created. By making a
finite number of such perturbations we can eliminate all successive
$5$-tangencies.

The same argument applies to the induction step.

Notice, that at each step a finite number of perturbations of three types is
made: the ones used to get that $3$-tangencies are smooth $1$-parameter
families of trajectories, the ones which eliminate non-transversal
intersections of ruled surfaces with the scatterers and the ones which
eliminate $5$-tangencies. One can choose the size of the perturbations at each
step in such a way that the perturbations made after step $N$ keep the
transversal manifolds considered at the previous steps transversal.

Indeed, let $\alpha_n$ be such that any perturbation of size less than
$\alpha_n$ keeps $3$-tangencies being smooth $1$-parameter families of
trajectories, $\beta_n<\alpha_n/2K^{n+5}$ be such that any perturbation of size
less than $\beta_n$ keeps ruled surfaces of order $\leq n$ transversal and let
$\gamma_n<\beta_n/2 K^{n+3}$ be such that any perturbation of size less than
$\gamma_k$ doesn't introduce $5$-tangencies of order $\leq k.$ One can choose
$\alpha_j,$ $\beta_j,$ $\gamma_j,$ $j=1,2,\ldots$ in such a way that for every
$N>0:$

$$\sum\limits_{n\geq N} K^{n+4}\alpha_n+K^{n+5}\beta_n+K^{n+3}\gamma_n\leq\gamma_{N-1}/2$$
and
$$\sum\limits_{n\geq 0} K^{n+4}\alpha_n+K^{n+5}\beta_n +K^{n+3}\gamma_n < \eps_0.$$

Then the series of perturbations converges, the total perturbations is less
that $\eps_0$ and there are no trajectories with more that $4$ tangencies for
the obtained collection of scatterers.
\end{proof}

From now on assume that for the configuration of the scatterers we consider
there are no trajectories with more that $2d-2$ tangencies.

\section{Estimates on the measure}

Let $B(x,r)$ denote the open ball of radius $r$ about $x$ in $\R^n.$ The
following theorem is proved in \cite{BF}.

\begin{theorem}
Let $F$ be a real-valued $C^m$ function on $B(0,1),$ with

\begin{enumerate}
     \item $c_0<\max\limits_{|\alpha|=m-1} |\da F(0)|< C_0,$ and with
     \item $|\da F|\leq C_1$ on $B(0,1)$ for $|\alpha|=m.$
        \label{hyptwothmone}
\end{enumerate}
 Let
\begin{enumerate}
    \item[(3)] $V(F)=\{x\in B(0,1): F(x)=0\},$ and let
    \item[(4)] $V(F,\delta)=\{x\in B(0,c_1):\
        {\rm distance}(x,V(F))<\delta\},$
\end{enumerate}
where $c_1$ is a small enough constant determined by $c_0,$ $C_0,$ $C_1,$ $m,$
$n.$

Then
\begin{equation}\nonumber
    Vol\{V(F,\delta)\}\leq C_2\delta\ {\rm for}\ 0<\delta<c_1,
\end{equation}
where $C_2$ is a large constant determined by $c_0,$ $C_0,$ $C_1,$ $m,$ $n.$

\label{measure}\end{theorem}

Thus, if $F$ is a smooth function that never vanishes to infinite order, then
the set of points within the distance $\delta$ from the zeroes of $F$ has
volume $O(\delta).$

\begin{corollary} Let $x_0\in\partial M\backslash\mathcal{S}.$ There exists a
neighborhood $U(x_0)\subset\dM\backslash\mathcal{S},$ such that the measure of
the $\delta$-neighborhood of any component of $T^{-k}\tang\cap U$ is of order
$O(\delta).$ \label{measure2}
\end{corollary}

\begin{proof}
 By Lemma \ref{tolines} there is a $C^{8d-7}$ diffeomorphism of a
neighborhood of $x_0$ to the neighborhood of the corresponding line in the
space of lines in $\R^d$ with the jacobian bounded away from zero.

From Lemma \ref{change}, Lemma \ref{resolution} and Proposition \ref{4ref} it
follows that in the neighborhood of any of its points in the space of lines the
singularity manifold $T^{-k}\tang$ is given by a $C^{4d-3}$ function with a
non-zero $4d-4$-jet. The theorem \ref{measure} gives the estimate in a
neighborhood of any point of $T^{-k}\tang\cap U.$
\end{proof}

\section{The fundamental theorem}
In this section we formulate a version of the fundamental theorem for
multidimensional dispersing billiards. This theorem was first formulated and
proved in \cite{S1} for plane billiards. It allows to apply the modified Hopf
argument to prove local ergodicity of the system. The derivation of the local
ergodicity from the fundamental theorem is done in details in \cite{KSSz},
\cite{LW} and is omitted here.

In the second part of this section we give the definitions required in the
proof of the fundamental theorem and show that under our main assumption in the
case of strictly dispersing scatterers one needs only finitely smooth
scatterers. The fundamental theorem for multidimensional dispersing billiards
was first formulated and proved by Chernov and Sinai \cite{SCh} and generalized
in \cite{KSSz}. We follow here the exposition in \cite{BCST2}.

\subsection{The statement of the Fundamental Theorem.}

Let $x\in \dM$ be such a point, that the trajectory of $x$ has at most one
tangency. Consider a small neighborhood $U(x).$ A parameterized family of
finite coverings $$\G=\{\Gi | i=1,\ldots, I(\delta)\}\qquad 0<\delta<\delta_0$$
is a family of regular coverings iff:
\begin{enumerate}
\item each $\Gi$ is a "parallelepiped", that is the image of a $2d-2$ dimensional
cube under the linear map $\R^{2d-2}\to \dM,$
\item the faces of each $\Gi$ are parallel to tangent planes of stable and unstable
manifolds  $W^{s}(w(\Gi))$ and $W^{u}(w(\Gi))$ passing through the center of
the parallelepiped $(w(\Gi)),$
\item For any point, the number of parallelepipeds covering it is at most $2^{2d-2},$
\item if $\Gi\cap G_j^{\delta}\neq\emptyset,$ then $$\mu(\Gi\cap
G_j^{\delta})\geq c_1 \delta^{2d-2}$$ where $c_1$ is independent of $\delta.$
\end{enumerate}

For a given $\Gi$ $s$-faces are those parallel to $W^{s}(w(\Gi))$ and $u$-faces
are those parallel to $W^{u}(w(\Gi)).$ For $y\in\Gi$ we say that a stable
manifold $W^s(y)$ intersects $\Gi$ {\it correctly}, if $\partial (\Gi\cap
W^{s}(y))$ belongs to the union of $u$-faces and that an unstable manifold
$W^u(y)$ intersects $\Gi$ {\it correctly}, if $\partial (\Gi\cap W^{u}(y))$
belongs to the union of $u$-faces.

\begin{theorem}
Let $x\in M$ be such a point, that the trajectory of $x$ has at most one
tangency and $0<\epso<1$ be a fixed constant. Then there exists a sufficiently
small neighborhood $U_{\epso},$ such that for any family of regular coverings
$\G$ of $U_{\epso},$ the covering $\G$ can be divided into two disjoint
subsets, $\G_g$ and $\G_b$ (called 'good' and 'bad'), in such a way that:
\begin{itemize}
\item For any $G_i^{\delta}\in\G_b$ and any $s$-face $E^s$ of it, the set
$$\{y\in G_i^{\delta} |\rho(y,E^s)<\epso\delta\ {\rm and}\ \gamma^s(y)\ {\rm
intersects}\ {\rm correctly}\}$$ has positive relative $\mu_1$ measure in
$G_i^{\delta}.$
\item $$\mu_1(\cup_{G_i^{\delta}\in\G_b} G_i^{\delta})=o(\delta)$$
\end{itemize}
\label{fundtheorem}
\end{theorem}

\subsection{The proof of the Fundamental theorem.}

Let us remind the main ideas of the proof of the theorem. Informally, the
parallelepiped is good if sufficiently many (in the sense of the relative
measure in the parallelepiped) stable manifolds are not cut inside this
parallelepiped. A given parallelepiped can be bad either because it intersects
the preimages of the singularity manifolds $\sing,\ttang,\ldots
T^{-F(\delta)}\mathcal{S}$ or because it intersects $T^{-F(\delta)-1},\ldots$
The set of bad parallelepipeds of the second type is small because of
hyperbolicity of the system. This is a so-called {\it tail estimate}, which
becomes rather simple in the setting of strictly dispersing scatterers (see
\cite{SCh} for the proof). As it was noticed in \cite{BCST2} this part doesn't
require information about the structure of the preimages of the singularity
manifolds. To estimate the measure of the remaining part of bad parallelepipeds
we notice that
\begin{enumerate}
    \item for a sufficiently small $\delta$ there are no parallelepipeds which intersect
        more than four singularity manifold from $\sing,\ttang,\ldots T^{-F(\delta)}\mathcal{S}.$
    \label{note1}
    \item a parallelepiped intersected by at most four singularity manifolds cannot be
        bad because by choosing a sufficiently small neighborhood of $x$ we can
        assume, that the angle between stable manifolds and singularity manifolds
        is arbitrary small.
    \label{note2}
\end{enumerate}

As it is explained in \cite{BCST2} it is here that the additional information
about singularity manifolds is needed.



We formulate the lemma about parallelization (see \cite{KSSz}, Lemma 4.9).

\begin{lemma} Given any $x\in \partial M^0$ and any $\eps>0$ there is
    a neighborhood $U(x)\subset\partial M$ such that for every
    $\gamma_1^s,\gamma_2^s$ and any 3-dimensional manifold $T^{-n}\tang$
    ($n>0$) intersecting $U(x)$ with points $y_1, y_2$ and $\hat{y},$ lying on
    three manifolds, respectively, so that $T_{\hat{y}}T^{-n}\tang$ exists:
    $$\angle(T_{y_1}\gamma_1^s,T_{y_2}\gamma_2^s)<\eps$$
    $$\angle(T_{y_1}\gamma_1^s,T_{\hat{y}}T^{-n}\tang)<\eps$$
\label{parallel}
\end{lemma}

We introduce some more notations. The following two quantities measure the
hyperbolicity near the point $y\in\partial M^0.$ Let
$$\kappa_{n,0}(y)=\inf\limits_{\Sigma} ||(D^n_{-T^n y,\Sigma})^{-1}||^{-1}_p,$$
where the infinum is taken over all convex local orthogonal manifolds passing
through $-T^n y.$ Furthermore denote
$$\kappa_{n,\delta}(y)=
    \inf\limits_{\Sigma}\inf\limits_{w\in\Sigma}||(D^n{w,\Sigma})^{-1}||^{-1}_p$$
Here the infinum is taken for the set of convex fronts $\Sigma$ passing through
$-T^n y$ such that $T^n$ is continuous on $\Sigma$ and $T^n\Sigma\subset
B_{\delta}(-y).$

The following subsets of the neighborhood $U(x)$ depend on the constant
$\delta:$
$$ U^g=\{y\in U | \forall n\in\Z_{+}, z(T^n
    y)\geq(\kappa_{n,c_3\delta}(y))^{-1} c_3\delta\},$$
$$U^b=U\backslash U^g$$
$$ U^b_n=\{y\in U | z_{tub}(T^n y) < (\kappa_{n,c_3\delta}(y))^{-1} c_3\delta\}$$





\begin{definition} A function $F:\R_{+}\to\Z_{+}$ defined in a neighborhood of
the origin is called permitted if $F(\delta)\uparrow\infty$ as
$\delta\downarrow 0.$ For a fixed permitted function $F(\delta)$ we define
$U^b_\omega=\cup_{n>\delta} U^b_n.$
\end{definition}

\begin{lemma}[Tail bound] For any permitted function $F(\delta):$
$$\mu_1(U^b_{\omega})=\bar{o}(\delta)$$
\label{tail}
\end{lemma}

The proof of the fundamental theorem follows from the following estimates:
\begin{itemize}

\item The estimate from below on the measure of $G^{\delta}_i\cap
    (E^s)^{[\eps_1\delta]},$ where  $E^s$ is an $s$-face for the bad parallelepiped
    $G_i^{\delta}:$
    \be
        \mu_1(G^{\delta}_i\cap (E^s)^{[\eps_1\delta]})\geq
        c\eps_1^{d-1}\mu_1(G^{\delta}_i)\geq\eps_3\mu_1 (G^{\delta}_i)
    \ee
    in case $\eps_3=\eps_3(\eps_1)$ is chosen sufficiently small.

\item Similarly,
    \be
        \mu_1(G_i^{\delta}\cap(\partial\Gi)^{[\eps_4\delta]})\leq
        \frac{\eps_3}{4}\mu_1(\Gi)
    \ee
for $\eps_4=\eps_4(\eps_3)$ chosen sufficiently small.

\item The estimate on the measure of the $\eps_4 \delta$-neighborhood of the
singularity manifold inside the parallelepiped. The Corollary \ref{measure2}
implies that \be
    \mu_1(\Gi\cap (\mathcal{S}^n)^{[\eps_4\delta]})\leq \eps_3/16\mu_1(\Gi)
\ee for $\eps_4=\eps_4(\eps_3)$ chosen sufficiently small.

\item We now choose $\eps_2(\eps_4)$ small enough, so that by Lemma \ref{parallel}
the stable manifolds and singularity components are 'almost parallel'.Namely,
the smallness of $\eps_2$ should guarantee that for any $y\in\Gi$ for which
$\gamma^s(y)$ doesn't intersect correctly we have:

\be
    y\in(\Gi\cap(\mathcal{S}^n)^{[\eps_4\delta]})\cup
    (G_i^{\delta}\cap(\partial\Gi)^{[\eps_4\delta]})\cup U^b_{\omega}
\ee

\end{itemize}

Let us now consider a bad parallelepiped $\Gi$ with an s-face $E^s$ for which

\be
    \mu_1(\Gi\cap(E^s)^{[\eps_1\delta]}\cap U_{ic})\leq\frac{\eps_3}{4}\mu_1(\Gi)
\label{bad} \ee

Here $U_{ic}$ is the set of points in $\Gi$ with correctly intersecting local
stable manifolds.

By the estimates above

$$\mu_1(\Gi\cap U^b_{\omega})\geq\mu_1(\Gi\cap(E^s)^{[\eps_1\delta]}\cap
U^b_{\omega}\geq\frac{\eps_3}{4}\mu_1(\Gi).$$

Now recall that in a regular covering there are at most $2^{2d-2}$
parallelepipeds with a non-empty common intersection. Thus:

$$2^{2d-2}\mu_1(U^b_{\omega})\geq\sum' \mu_1(\Gi\cap U^b_{\omega})\geq
    \frac{\eps_3}{4}\sum' \mu_1(\Gi),$$

where $\sum '$ denotes the sum over bad parallelepipeds for which (\ref{bad})
holds for some s-face. By the Tail Bound (Lemma \ref{tail}) we have $\sum '
\mu_1(\Gi)=\bar{o}(\delta)$ thus the proof of the theorem \ref{fundtheorem} is
complete.

\end{document}